\documentclass{amsart}
\usepackage{amssymb,amsmath}

\newtheorem{theorem}{Theorem}
\newtheorem{lemma}{Lemma}
\newtheorem{example}{Example}
\newtheorem{problem}{Problem}
\newcommand{\bt}{\begin{theorem}}
\newcommand{\et}{\end{theorem}}
\newcommand{\bl}{\begin{lemma}}
\newcommand{\el}{\end{lemma}}
\newcommand{\bex}{\begin{example}}
\newcommand{\eex}{\end{example}}
\newcommand{\bp}{\begin{problem}}
\newcommand{\ep}{\end{problem}}
\newcommand{\beal}{\begin{align*}}
\newcommand{\enal}{\end{align*}}
\newcommand{\beq}{\begin{equation}}
\newcommand{\eeq}{\end{equation}}
\newcommand{\benum}{\begin{enumerate}}
\newcommand{\eenum}{\end{enumerate}}
\newcommand{\ba}{\begin{array}}
\newcommand{\ea}{\end{array}}

\newcommand{\N}{\ensuremath{\mathbf N}}

\newcommand{\mcf}{\ensuremath{\mathcal F}}

\newcommand{\mcr}{\ensuremath{\mathcal R}}
\newcommand{\mcs}{\ensuremath{\mathcal S}}
\newcommand{\mct}{\ensuremath{\mathcal T}}
\newcommand{\mcu}{\ensuremath{\mathcal U}}
\newcommand{\mcv}{\ensuremath{\mathcal V}}
\newcommand{\mcw}{\ensuremath{\mathcal W}}
\newcommand{\mcrp}{\ensuremath{\mathcal R'}}
\newcommand{\mcsp}{\ensuremath{\mathcal S'}}
\newcommand{\mctp}{\ensuremath{\mathcal T'}}
\newcommand{\mcup}{\ensuremath{\mathcal U'}}
\newcommand{\mcvp}{\ensuremath{\mathcal V'}}
\newcommand{\mcwp}{\ensuremath{\mathcal W'}}
\newcommand{\mcrpp}{\ensuremath{\mathcal R''}}
\newcommand{\mcspp}{\ensuremath{\mathcal S''}}
\newcommand{\mctpp}{\ensuremath{\mathcal T''}}

\newcommand{\polf}{$\mathcal{F} = \{f_n(q)\}_{n=1}^{\infty}$}

\newcommand{\polrp}{$\mathcal{R}' = \{r'_n(q)\}_{n=1}^{\infty}$}
\newcommand{\polsp}{$\mathcal{S}' = \{s'_n(q)\}_{n=1}^{\infty}$}
\newcommand{\poltp}{$\mathcal{T}' = \{t'_{m,n}(q)\}_{m,n=1}^{\infty}$}
\newcommand{\polu}{$\mathcal{U} = \{u_n(q)\}_{n=1}^{\infty}$}
\newcommand{\polv}{$\mathcal{V} = \{v_m(q)\}_{m=1}^{\infty}$}

\newcommand{\polup}{$\mathcal{U}' = \{u'_n(q)\}_{n=1}^{\infty}$}
\newcommand{\polvp}{$\mathcal{V}' = \{v'_m(q)\}_{m=1}^{\infty}$}
\newcommand{\polwp}{$\mathcal{W}' = \{w'_{m,n}(q)\}_{m,n=1}^{\infty}$}

\begin{document}

\title{Quadratic addition rules for quantum integers}
\subjclass[2000]{Primary 30B12, 81R50.  Secondary 11B13.}
\keywords{Quantum integers, quantum polynomial, quantum addition,
polynomial functional equation, $q$-series.}
\author{Alex V. Kontorovich}
\address{Department of Mathematics\\
Columbia University\\
New York, New York 10027}
\email{alexk@math.columbia.edu}
\author{Melvyn B. Nathanson}
\address{Department of Mathematics\\
Lehman College (CUNY)\\
Bronx, New York 10468}
\email{melvyn.nathanson@lehman.cuny.edu}
\thanks{The work of M.B.N. was supported in part by grants from the NSA
Mathematical Sciences Program and the PSC-CUNY Research Award
Program.}

\begin{abstract}
For every positive integer $n$, the quantum integer $[n]_q$ is the
polynomial $[n]_q = 1 + q + q^2 + \cdots + q^{n-1}.$ A quadratic
addition rule for quantum integers consists of  sequences of
polynomials $\mathcal{R}' = \{r'_n(q)\}_{n=1}^{\infty}$,
$\mathcal{S}' = \{s'_n(q)\}_{n=1}^{\infty}$, and $\mathcal{T}' =
\{t'_{m,n}(q)\}_{m,n=1}^{\infty}$ such that $[m+n]_q = r'_n(q)[m]_q
+ s'_m(q)[n]_q + t'_{m,n}(q)[m]_q[n]_q$ for all $m$ and $n.$  This
paper gives a complete classification of quadratic addition rules,
and also considers sequences of polynomials \polf\ that satisfy the
associated functional equation $f_{m+n}(q)= r'_n(q)f_m(q) +
s'_m(q)f_n(q) + t'_{m,n}f_m(q)f_n(q).$
\end{abstract}

\maketitle

\section{Quantum addition rules}
Let \N\ denote the set of positive integers.
Let $K[q]$ denote the ring of polynomials with coefficients in a field $K$.
For every positive integer $n$, the quantum integer $[n]_q$ is the polynomial
\[
[n]_q = 1 + q + q^2 + \cdots + q^{n-1} \in K[q].
\]
We define $[0]_q = 0  .$
Equivalently,
\[
[n]_q  = \frac{1-q^n}{1-q}.
\]
With ordinary addition of polynomials, $[m]_q + [n]_q \neq [m+n]_q.$
Our goal is to describe a class of binary operations $\oplus$ of ``quantum addition''
on the sequence $\{ [n]_q \}_{n=1}^{\infty}$ of quantum integers such that
\[
[m]_q \oplus [n]_q = [m+n]_q.
\]

A general setting for these binary operations is as follows.
Suppose that for every pair $(m,n)$ of positive integers there is a function
\[
\Phi_{m,n}:K[q]\times K[q] \rightarrow K[q].
\]
We use the sequence $\{\Phi_{m,n}\}_{m,n=1}^{\infty}$ to construct a binary operation $\oplus$ on the elements of an arbitrary sequence \polf\ of polynomials in $K[q]$.  For every pair $(m,n)$ of positive integers we define
\[
f_m(q) \oplus f_n(q) = \Phi_{m,n}(f_m(q),f_n(q)).
\]
The binary operation defined by $\{\Phi_{m,n}\}_{m,n=1}^{\infty}$ is called a {\em polynomial addition rule}.
A polynomial addition rule is {\em consistent} on the sequence \polf\ if
\[
f_m(q) \oplus f_n(q) = f_{r}(q)\oplus f_s(n)
\]
for all positive integers $m,n,r,$ and $s$ such that $m+n = r+s.$
If $\oplus$ is consistent on the sequence \polf\, then $\oplus$ is commutative in the sense that
\[
f_m(q) \oplus f_n(q) = f_{n}(q)\oplus f_m(n)
\]
for all positive integers $m$ and $n$.

A {\em quantum addition rule} is a polynomial addition rule that
is {\em consistent} on the sequence of quantum integers and
satisfies the functional equation
\beq          \label{Quad:afe}
[m]_q \oplus [n]_q = [m+n]_q
\eeq
for all positive integers $m$ and $n$.
For example, let
\beq           \label{Quad:exlinear}
\Phi_{m,n}(a(q),b(q)) = a(q) + q^m b(q)
\eeq
for all polynomials $a(q), b(q) \in K[q].$  This defines a quantum addition rule, since
\[
\begin{split}
\Phi_{m,n}([m]_q,[n]_q) & = [m]_q + q^m [n]_q  \\
& (1 + q + \cdots + q^{m-1}) + q^m(1 + q + \cdots + q^{n-1}) \\
& = 1 + q +\cdots + q^{m-1} + q^m +  \cdots + q^{m+n-1}  \\
& = [m+n]_q.
\end{split}
\]

There is also a functional equation associated to every quantum addition rule:
Find all sequences \polf\ of polynomials such that $\oplus$ is consistent on $\mathcal{F}$,
and $\mathcal{F}$ satisfies
\beq   \label{Quad:assocfe}
f_{m+n}(q) = f_m(q) \oplus f_n(q) = \Phi_{m,n}(f_m(q),f_n(q))
\eeq
for all positive integers $m$ and $n$.
For the quantum addition rule~(\ref{Quad:exlinear}),  the sequence of polynomials \polf\ satisfies the associated functional equation
\[
f_{m+n}(q) = f_m(q) + q^m f_n(q)
\]
if and only if there exists a polynomial $h(q)$ such that
\[
f_n(q) = [n]_qh(q)
\]
for all $n$ (Nathanson~\cite{nath05a}).  The analogous functional equation associated to multiplication of quantum integers has been studied by Nathanson~\cite{nath03b,nath04c} and Borisov, Nathanson, and Wang~\cite{bori-nath-wang04}.

A quantum addition rule is {\em linear} if there exist sequences of polynomials
\polrp\ and \polsp\ such that
\[
\Phi_{m,n}(a(q),b(q)) = r'_n(q)a(q) + s'_m(q)b(q)
\]
and so
\[
[m+n]_q = r'_n(q)[m]_q + r'_m(q)[n]_q
\]
for all $m,n \in \N.$
The quantum addition rule~(\ref{Quad:exlinear}) is linear.
Linear quantum addition rules are considered in Nathanson~\cite{nath05a,nath05b}.

A quantum addition rule is {\em quadratic} if there exist sequences of polynomials
\polrp, \polsp, and \poltp\ such that
\[
[m+n]_q = r'_n(q)[m]_q + s'_m(q)[n]_q + t'_{m,n}(q)[m]_q[n]_q
\]
for all $m,n \in \N.$
Examples of quadratic quantum addition rules are
\beq   \label{Quad:rule1}
[m+n]_q =  [m]_q + [n]_q - (1 - q)[m]_q[n]_q
\eeq
and
\beq     \label{Quad:rule2}
[m+n]_q = q^n [m]_q + q^m[n]_q + (1-q)[m]_q[n]_q.
\eeq
The associated functional equations are
\[
f_m(q)\oplus f_n(q) =  f_m(q) + f_n(q) - (1 - q)f_m(q)f_n(q)
\]
and
\[
f_m(q)\oplus f_n(q) = q^n f_m(q) + q^mf_n(q) + (1-q)f_m(q)f_n(q).
\]
In this paper we shall give a complete classification of quadratic addition rules for quantum integers, and discuss their associated functional equations.

\section{Quadratic zero identities}
We begin with zero identities.
A {\em quadratic zero identity} consists of three sequences of polynomials \polup, \polvp, and \polwp\ such that
\[
u'_n(q)[m]_q + v'_m(q)[n]_q + w'_{m,n}(q)[m]_q[n]_q = 0
\]
for all positive integers $m$ and $n$.
The following theorem classifies all quadratic zero identities.

\bt  \label{Quad:theorem:zero}
The sequences of polynomials  \polup, \polvp, and \polwp\ satisfy the quadratic zero identity
\beq      \label{Quad:zero}
u'_n(q)[m]_q + v'_m(q)[n]_q + w'_{m,n}(q)[m]_q[n]_q = 0
\eeq
for all positive integers $m$ and $n$ if and only if there exist sequences \polu\ and \polv\ such that
\beq        \label{Quad:u}
u'_n(q) = u_n(q)[n]_q,
\eeq
\beq        \label{Quad:v}
v'_m(q) = v_m(q)[m]_q,
\eeq
and
\beq        \label{Quad:w}
w'_{m,n}(q) = -(u_n(q)+v_m(q))
\eeq
for all positive integers $m$ and $n$.
\et

\begin{proof}
Suppose that the sequences \mcu', \mcv', and \mcw' satisfy~(\ref{Quad:zero}).
For all positive integers $m$ and $n$ we define
\[
u_n(q) = -( v'_1(q) + w'_{1,n}(q) )
\]
and
\[
v_m(q) = -( u'_1(q) + w'_{m,1}(q) ).
\]
Letting $m = 1$ in equation~(\ref{Quad:zero}), we obtain
\[
u'_n(q) + v'_1(q)[n]_q + w'_{1,n}(q)[n]_q = 0,
\]
and so
\[
u'_n(q) = -( v'_1(q) + w'_{1,n}(q) ) [n]_q = u_n(q)[n]_q.
\]
Letting $n = 1$ in equation~(\ref{Quad:zero}), we obtain
\[
u'_1(q)[m]_q + v'_m(q) + w'_{m,1}(q)[m]_q = 0,
\]
and so
\[
v'_m(q) = -( u'_1(q) + w'_{m,1}(q) ) [m]_q = v_m(q)[m]_q.
\]
Inserting these expressions for $u'_n(q)$ and $v'_m(q)$ into~(\ref{Quad:zero}) and dividing by $[m]_q[n]_q$, we obtain $w'_{m,n}(q) = -(u_n(q)+v_m(q).$  This proves that every zero identity is obtained from a pair of sequences \mcu\ and \mcv\ by the  construction~(\ref{Quad:u}),~(\ref{Quad:v}), and~(\ref{Quad:w}).

It is an immediate verification that we obtain a quadratic zero identity by applying~(\ref{Quad:u}),~(\ref{Quad:v}), and~(\ref{Quad:w}) to any two sequences \mcu\ and \mcv.
\end{proof}

A {\em linear zero identity} consists of two sequences of polynomials
\polup\ and \polvp  such that
\[
u'_n(q)[m]_q + v'_m(q)[n]_q = 0
\]
for all positive integers $m$ and $n$.

\bt
The sequences of polynomials \polup\ and \polvp satisfy the linear zero identity
\beq      \label{Quad:linearzero}
u'_n(q)[m]_q + v'_m(q)[n]_q = 0
\eeq
for all positive integers $m$ and $n$ if and only if there is a polynomial $z(q)$such that
\beq        \label{Quad:lin-u}
u'_n(q) = z(q)[n]_q
\eeq
and
\beq         \label{Quad:lin-v}
v'_m(q) = -z(q)[m]_q
\eeq
for all positive integers $m$ and $n$.
\et

\begin{proof}
The linear zero identity~(\ref{Quad:linearzero}) is a quadratic zero identity with $w'_{m,n}(q)=0$ for all $m$ and $n$.
It follows from Theorem~\ref{Quad:theorem:zero} that there exist polynomials $u_n(q)$ and $v_m(q)$ such that
\[
w'_{m,n}(q) = -(u_n(q)+v_m(q)) = 0,
\]
and so there exists a polynomial $z(q)$ such that
\[
u_n(q) = -v_m(q) = z(q),
\]
\[
u'_n(q) = z(q)[n]_q
\]
and
\[
v'_m(q) = -z(q)[m]_q
\]
for all $m$ and $n$.
Conversely, if there exists a polynomial $z(q)$ such that the sequences \polup\ and \polvp\ satisfy~(\ref{Quad:lin-u}) and~(\ref{Quad:lin-v}), then we obtain the zero identity~(\ref{Quad:linearzero}).
\end{proof}

The following result follows immediately from Theorem~\ref{Quad:theorem:zero}.

\bt
Let \polup, \polvp, and \polwp\ be sequences of polynomials that satisfy the zero identity
\[
u'_n(q)[m]_q + v'_m(q)[n]_q + w'_{m,n}(q)[m]_q[n]_q = 0
\]
for all positive integers $m$ and $n$.
If $\deg(u'_n(q)) < n-1,$ then $u'_n(q) = 0.$
If $\deg(v'_m(q)) < m-1,$ then $v'_m(q) = 0.$
\et

\section{Quadratic addition rules}
A {\em quadratic addition rule} for the quantum integers consists of three sequences of polynomials
\polrp, \polsp, and \poltp\ such that
\beq                                          \label{Quad:addrule}
[m+n]_q = r'_n(q)[m]_q + s'_m(q)[n]_q + t'_{m,n}(q)[m]_q[n]_q
\eeq
for all positive integers $m$ and $n$.
Our goal is to classify all quadratic addition rules, that is, to find all sequences of polynomials $\mathcal{R}',$ $\mathcal{S}',$ and $\mathcal{T}'$
that satisfy~(\ref{Quad:addrule}) for all quantum integers $[m]_q$ and $[n]_q.$

Suppose that the sequences \mcrp, \mcsp, and \mctp\ determine a quadratic addition rule,
and that the sequences \mcrpp, \mcspp, and \mctpp\ also determine a quadratic addition rule.
Then
\[
(r'_n(q)-r''_n(q))[m]_q + (s'_m(q)-s''_m(q))[n]_q + (t'_{m,n}(q)-t''_{m,n}(q))=0,
\]
and so the sequences \mcup = \mcrp - \mcrpp, \mcvp = \mcsp - \mcspp, and \mcwp = \mctp - \mctpp\ determine a quadratic zero identity.
Similarly, if the sequences \mcup, \mcvp, and \mcwp\ determine a quadratic zero identity, then for every scalar $\lambda$ the sequences $\mcrp+\lambda\mcup$, $\mcsp+\lambda\mcvp$, and $\mctp + \lambda\mcwp$ also determine a quadratic addition rule.
Thus, every quadratic addition rule for the quantum integers can be expressed as the sum of a fixed rule and a zero identity.

We can use polynomial division to find a standard form for a quadratic addition rule.  Let \mcrp, \mcsp, and \mctp\ be sequences of polynomials that satisfy~(\ref{Quad:addrule}).  By the division algorithm for polynomials, for every $n \geq 1$ there exist polynomials $u_n(q)$ and $r_n(q)$ such that
\[
r'_n(q) = u_n(q)[n]_q + r_n(q)
\]
and
\[
\deg(r_n(q)) \leq \deg([n]_q) - 1 = n-2.
\]
Similarly, for every $m \geq 1$ there exist polynomials $v_m(q)$ and $s_m(q)$ such that
\[
s'_m(q) = v_m(q)[m]_q + s_m(q)
\]
and
\[
\deg(s_m(q)) \leq \deg([m]_q) - 1 =  m-2.
\]
We define the polynomials
\[
u'_n(q) = u_n(q)[n]_q,
\]
\[
v'_m(q) = v_m(q)[m]_q,
\]
and
\[
w'_{m,n}(q) = -(u_n(q)+v_m(q)),
\]
and obtain the quadratic zero identity
\beq        \label{Quad:zerorule2}
0 = u'_n(q)[m]_q + v'_m(q)[n]_q + w'_{m,n}(q)[m]_q[n]_q.
\eeq

Let
\[
t_{m,n}(q) = t'_{m,n}(q) - w'_{m,n}(q) = t'_{m,n}(q) + u_n(q)+v_m(q).
\]
Subtracting the quadratic zero identity~(\ref{Quad:zerorule2}) from the quadratic addition rule~(\ref{Quad:addrule}), we obtain a new quadratic addition rule
\[
\begin{split}
[m+n]_q
& = (r'_n(q)-u'_n(q))[m]_q + (s'_m(q)-v'_m(q))[n]_q + (t'_{m,n}(q)-w'_{m,n}(q))[m]_q[n]_q \\
& = r_n(q)[m]_q + s_m(q)[n]_q + t_{m,n}(q)[m]_q[n]_q.
\end{split}
\]
The degrees of the polynomials in this identity satisfy
\[
\deg([m+n]_q = m + n-1,
\]
\[
\deg(r_n(q)[m]_q) \leq m+n-3,
\]
\[
\deg(s_m(q)[n]_q) \leq m+n-3,
\]
and
\[
\deg(t_{m,n}(q)[m]_q[n]_q) \geq m+n-2.
\]
This implies that
\[
\deg(t_{m,n}(q)) = 1
\]
Moreover, $t_{m,n}(q)$ is a monic polynomial since the quantum integers $[m]_q, [n]_q,$ and $[m+n]_q$ are monic.

For example, if we apply this procedure to the quadratic addition rules
\[
[m+n]_q = [m]_q + [n]_q + (q-1)[m]_q[n]_q
\]
or
\[
[m+n]_q = q^n[m]_q + q^m[n]_q - (q-1)[m]_q[n]_q,
\]
we obtain the rule
\[
[m+n]_q = r_n(q)[m]_q + s_m(q)[n]_q + t_{m,n}(q)[m]_q[n]_q,
\]
where
\[
r_n(q) = 1 - \delta_n,
\]
\[
s_m(q) = 1 - \delta_m,
\]
\[
t_{m,n}(q) = q-1 + \delta_m + \delta_n,
\]
and
\[
\delta_n = \left\{\ba{ll}
1 & \text{if $n = 1$,} \\
0 & \text{if $n \geq 2.$}
\ea
\right.
\]
This is called the {\em fundamental quadratic addition rule}.

Combining this with Theorem~\ref{Quad:theorem:zero}, we obtain the following complete classification of quadratic addition rules for quantum integers.

\bt      \label{Quad:theorem:quadrule}
Let \polrp, \polsp, and \poltp\ be sequences of polynomials.  Then
\[
[m+n]_q = r'_n(q)[m]_q + s'_m(q)[n]_q + t'_{m,n}(q)[m]_q[n]_q
\]
for all positive integers $m$ and $n$ if and only if there exist sequences of polynomials
\polu\ and \polv\ such that
\[
r'_n(q) = u_n(q)[n]_q + 1 - \delta_n,
\]
\[
s'_m(q) = v_m(q)[m]_q +1 - \delta_m,
\]
and
\[
t'_{m,n}(q) = q-1 - u_n(q)-v_m(q) +\delta_m + \delta_n.
\]
\et

\section{Functional equations associated to quadratic addition rules}
Let \mcr', \mcs', and \mct'\ define a quadratic rule for quantum addition, that is,  \polrp, \polsp, and \poltp\ are sequences of polynomials in $K[q]$ such that
\[
[m+n]_q = r'_n(q)[m]_q + s'_m(q)[n]_q + t'_{m,n}(q)[m]_q[n]_q
\]
for all positive integers $m$ and $n$.
Let \polf\ be a sequence of polynomials.  We define an addition operation $\oplus: \mcf \times \mcf \rightarrow K[q]$ by
\[
f_m(q) \oplus f_n(q) = r'_n(q)f_m(q) + s'_m(q)f_n(q) + t'_{m,n}(q)f_m(q)f_n(q).
\]
We want to find all sequences \mcf\ of polynomials such that, for all positive integers $m$ and $n$, we have
\[
f_{m+n}(q) = f_m(q) \oplus f_n(q)
\]
or, equivalently,
\beq           \label{Quad:addfe}
f_{m+n}(q) = r'_n(q)f_m(q) + s'_m(q)f_n(q) + t'_{m,n}(q)f_m(q)f_n(q).
\eeq
This functional equation always has the solution $f_n(q) = [n]_q$ for all $n\in \N,$ and also the zero solution $f_n(q) = 0$ for all $n \in \N.$
We would like to find every solution of this nonlinear equation.

Note that if the sequence of polynomials \polf\ is a solution of the functional equation~(\ref{Quad:addfe}), then \polf\ is determined inductively by the polynomial $f_1(q) = h(q)$, since
\beq        \label{Quad:induct}
\begin{split}
f_n(q) & = f_{1+(n-1)}(q) \\
& = h(q) \oplus f_{n-1}(q) \\
& = r'_{n-1}(q)h(q) + s'_1(q)f_{n-1}(q) + t'_{1,n-1}(q)h(q)f_{n-1}(q)
\end{split}
\eeq
for all $n \geq 2.$
Equivalently, we can ask for what polynomials $h(q)$ is the sequence \polf\ constructed from~(\ref{Quad:induct}) a solution of the functional equation~(\ref{Quad:addfe}).
We know that the polynomials $h(q) = 0$ and $h(q) = 1$ always produce solutions.

We can compute explicit solutions for the three quadratic addition rules discussed in this paper.  For each of these rules there is a solution of the functional equation with $f_1(q) = h(q)$ for every $h(q) \in K[q].$  Associated to the fundamental quadratic addition rule
\beq        \label{Quad:rule0}
[m+n]_q = (1-\delta_n)[m]_q + (1 - \delta_m)[n]_q + (q-1 + \delta_m + \delta_n)[m]_q[n]_q
\eeq
is the functional equation
\[
f_{m+n}(q) = (1-\delta_n)f_m(q) + (1 - \delta_m)f_n(q) + (q-1 + \delta_m + \delta_n)f_m(q)f_n(q).
\]
If \polf\ is a solution of this equation with $f_1(q) = h(q)$, then
\[
f_2(q) = (q+1)h(q)^2
\]
and
\[
\begin{split}
f_n(q) & = h(q) + q h(q)f_{n-1}(q) \\
& = h(q) + qh(q)^2 + q^2h(q)^3 + \cdots + q^{n-3}h(q)^{n-2} + \left(q^{n-2}+q^{n-1}\right)h(q)^n  \\
& = \frac{h(q)\left(1 - (qh(q))^{n} \right)}{1-qh(q)} + q^{n-2}h(q)^{n-1}\left(h(q)-1\right)
\end{split}
\]
for all $n \geq 3.$

For the quadratic addition rule~(\ref{Quad:rule1})
\[
[m+n]_q = [m]_q + [n]_q - (1-q)[m]_q[n]_q,
\]
the general solution of the functional equation
\[
f_{m+n}(q) = f_m(q) + f_n(q) - (1-q)f_m(q)f_n(q)
\]
is
\[
f_n(q) = \frac{1-(1-(1-q)h(q))^{n}}{1-q}.
\]
For the quadratic addition rule~(\ref{Quad:rule2})
\[
[m+n]_q = q^n[m]_q + q^m[n]_q + (1-q)[m]_q[n]_q,
\]
the general solution of the functional equation
\[
f_{m+n}(q) = q^nf_m(q) + q^mf_n(q) + (1-q)f_m(q)f_n(q)
\]
is
\[
f_n(q) = \frac{ (q+(1-q)h(q))^n - q^n }{1-q}.
\]

By Theorem~\ref{Quad:theorem:quadrule}, the general functional equation associated to a quadratic quantum addition rule is
\beq \label{Quad:genfe}
\begin{split}
f_{m+n}(q) = & \left( u_n(q)[n]_q + 1 - \delta_n \right) f_m(q) + \left(  v_m(q)[m]_q +1 - \delta_m\right)f_n(q) \\
& + \left( q-1 - u_n(q)-v_m(q) +\delta_m + \delta_n\right)f_m(q)f_n(q),
\end{split}
\eeq
for all positive integers $m$ and $n$, where \polu\ and \polv\ are arbitrary sequences of polynomials in $K[q].$

Let \polf\ be a solution of~(\ref{Quad:genfe}).  This sequence is generated by the polynomial $f_1(q) = h(q).$  Setting $m = n = 1$, we obtain
\beq    \label{Quad:f2}
f_2(q) = (u_1(q)+v_1(q))h(q) + (q+1-u_1(q)-v_1(q))h(q)^2.
\eeq
Setting $m = 1$ and $n=2$, we obtain
\[
\begin{split}
f_3(q)
& = f_{1+2}(q) \\
& = \left( u_2(q)(q+1) + 1 \right) h(q) + v_1(q)f_2(q)
+ \left( q - u_2(q)-v_1(q) \right)h(q)f_2(q).
\end{split}
\]
Similarly, with $m=2$ and $n=1$, we obtain
\[
\begin{split}
f_3(q) & = f_{2+1}(q) \\
& = u_1(q)f_2(q) + \left(  v_2(q)(q+1) +1 \right)h(q)
+ \left( q - u_1(q)-v_2(q) \right)f_2(q)h(q),
\end{split}
\]
Subtracting these equations gives
\[
\begin{split}
0 = &
\left(  (u_1(q)-v_1(q)) - (u_2(q)-v_2(q)) \right) f_2(q)h(q) +\\
& + \left( u_2(q)-v_2(q)  \right)(q+1) h(q) + (v_1(q) - u_1(q)) f_2(q).
\end{split}
\]
Replacing $f_2(q)$ by~(\ref{Quad:f2}), we see that the polynomial $h(q)$ must satisfy the identity
\[
\begin{split}
0 = & (u_2-v_2-u_1+v_1)(q+1-u_1-v_1)h(q)^3 \\
& + ((u_2-v_2)(u_1+v_1)+(q+1)(u_1-v_1)-2(u_1^2-v_1^2))h(q)^2 \\
& + (u_1^2-v_1^2-(u_2-v_2)(q+1))h(q).
\end{split}
\]
Equivalently, $h(q)$ is a root of the cubic polynomial
\[
\begin{split}
0 = & (u_2-v_2-u_1+v_1)(q+1-u_1-v_1)x^3 \\
& + ((u_2-v_2)(u_1+v_1)+(q+1)(u_1-v_1)-2(u_1^2-v_1^2))x^2 \\
& + (u_1^2-v_1^2-(u_2-v_2)(q+1))x.
\end{split}
\]
We know that $x = 0$ and $x = 1$ are solutions of this equation.  Dividing by $x(x-1),$ we obtain
\[
\begin{split}
(u_2-v_2-u_1+v_1)(q+1-u_1-v_1)x - (u_1^2-v_1^2-(u_2-v_2)(q+1)) = 0.
\end{split}
\]
For the quadratic addition rules~(\ref{Quad:rule0}),~(\ref{Quad:rule1}), and~(\ref{Quad:rule2}), the coefficients in this equation are both 0, and every polynomial $h(q)$ is a solution.
In general, however, if the coefficient of $x$ is nonzero, then this equation has at most one solution $h(q) \in K[q].$

It is an open problem to determine all solutions of the functional equations associated to quadratic zero rules and quadratic zero identities.
It is also of interest to find solutions of these functional equations in the ring of formal power series $K[[q]].$

\providecommand{\bysame}{\leavevmode\hbox to3em{\hrulefill}\thinspace}
\providecommand{\MR}{\relax\ifhmode\unskip\space\fi MR }
\providecommand{\MRhref}[2]{%
  \href{http://www.ams.org/mathscinet-getitem?mr=#1}{#2}
}
\providecommand{\href}[2]{#2}

\end{document}